\documentclass[draft]{article}
\usepackage{amsmath,amsthm}
 \usepackage{amssymb}
 \usepackage{amsfonts}
 \usepackage{mathrsfs}
 \usepackage{amssymb,amsmath,graphicx,graphics,amsthm,amsfonts,mathrsfs,multicol,multirow, amscd,color}
 \usepackage{graphicx}

\chardef\bslash=`\\

\newcommand {\red} {\textcolor{red}}
\newcommand {\green} {\textcolor{green}}
\newcommand {\blue} {\textcolor{blue}}

\newtheorem{thm}{Theorem}[section]
\newtheorem{cor}[thm]{Corollary}
\newtheorem{lem}[thm]{Lemma}
\newtheorem{pro}[thm]{Proposition}
\newtheorem{prob}[thm]{Problem}

\theoremstyle{definition}

\theoremstyle{remark}

\newcommand{\X}{\mathcal{X}}

\newcommand{\eval}[2][\right]{\relax
  \ifx#1\right\relax \left.\fi#2#1\rvert}

\def\proof{{\bf Proof.}\ }
\def\endproof{\hfill \ \ \hbox{$\sqcup$}\llap{\hbox{$\sqcap$}}}
\def \usecolour
{\newcommand {\rered}{\red}
\newcommand {\reblue} {\blue}
\newcommand {\regreen} {\green}
}

\usecolour

\newcommand {\relabel}[1] {\label{#1} \red{[*: #1]}}
\newcommand {\rebibitem}[1] {\bibitem{#1} \red{[*: #1]}}
\def\relabel {\label} \def\rebibitem {\bibitem}

\def \modtwo {{\ (\mbox {mod}\ 2)}}
\def \C {{\cal C}}

\def \C {{\cal C}}
\def \E {{\cal E}}

\def \E {{\cal E}}
\def \L {{\cal L}}
\def \M {{\cal M}}

\def \X {{\cal X}}

\def \ssim {\stackrel s\thicksim}
\def \nssim {\stackrel s\nsim}

\def \iff {if and only if }

\def \F {{\cal F}}
\def \nF {\bar {\cal F}}
\def \modtwo {\ (\mbox{mod } 2)}

\newcommand \Union[2]
{
\mbox{Union}_{#1}(#2)
}

\begin{document}
\title{On non-feasible edge sets in matching-covered graphs}
\author{Xiao Zhao\thanks{Corresponding author. Department of Mathematics, Harbin Institute of Technology, Harbin 150001, China. Email: zhaoxiao05@126.com},
Fengming Dong\footnote{National Institute  of Education,
Nanyang Technological University,
Singapore. Email:  fengming.dong@nie.edu.sg},
Sheng Chen\footnote{
Department of Mathematics, Harbin Institute of Technology, Harbin 150001, China. Email:  schen@hit.edu.cn}}
\date{November 28, 2018}

\maketitle
\markboth{On non-feasible edge sets in matching-covered graphs}
{}
\renewcommand{\sectionmark}[1]{}

\begin{abstract}
Let $G=(V,E)$ be a matching-covered graph and
$X$ be an edge set of $G$.
$X$ is said to be feasible if there
exist two perfect matchings $M_1$ and $M_2$ in $G$
such that $|M_1\cap X|\not \equiv|M_2\cap X|\ (\mbox{mod } 2)$.
For any $V_0\subseteq V$, $X$ is said to be
switching-equivalent to $X\oplus \nabla_G(V_0)$,
where $\nabla_G(V_0)$ is the set of edges in $G$
each of which has exactly one end in $V_0$
and $A \oplus B$ is
the symmetric difference of two sets $A$ and $B$.
Lukot'ka and Rollov\'a
showed that
when $G$ is  regular and bipartite,
$X$ is non-feasible if and only if
$X$ is switching-equivalent to $\emptyset$.
This article extends Lukot'ka and Rollov\'a's result by showing that this conclusion
holds as long as $G$ is matching-covered and bipartite.
This article also studies matching-covered graphs $G$
whose non-feasible edge sets are switching-equivalent to
$\emptyset$ or $E$ and partially characterizes these matching-covered
graphs in terms of their ear decompositions.
Another aim of this article is to construct infinite many
$r$-connected and $r$-regular graphs of class 1
containing non-feasible edge sets not switching-equivalent to either $\emptyset$ or $E$ for an arbitrary integer $r$ with $r\ge 3$,
which provides negative answers to problems
asked by Lukot'ka and Rollov\'a
and He, et al
respectively.
\end{abstract}

\section{Introduction and Preliminary}

This article studies
finite and undirected loopless graphs.
Let $G=(V,E)$ be a graph.
A {\it perfect matching} of $G$ is a set of independent edges which covers all vertices of $G$.
$G$ is said to be  {\it matching-covered}
if it is connected and each edge of $G$ is contained in
some perfect matching of $G$.
It is not difficult to verify that
any regular graph of class 1 is matching-covered.

For a matching-covered graph $G$,
an edge set $X$ of $G$ is
said to be {\it feasible}
if $G$ has two perfect matchings $M_1$ and $M_2$ such that
$|M_1\cap X|\not\equiv|M_2\cap X|\ (\mbox{mod } 2)$ holds.
Thus an edge set $X$ of $G$ is non-feasible
\iff $|M_1\cap X|\equiv|M_2\cap X|\ (\mbox{mod } 2)$
holds for every pair of perfect matchings
$M_1$ and $M_2$ of $G$.
For example, $E$ and $\emptyset$ are non-feasible edge sets
of $G$.
In Theorem~\ref{main2},
we extend the definition of a feasible edge
to connected graphs which are not matching-covered.

For any $V_0\subseteq V$,
let $\nabla_G(V_0)$ be the set of edges in $G$
each of which has exactly one end in $V_0$.
For any vertex $v$ in $G$,
$\nabla_G(\{v\})$ is exactly the set of edges in $G$
which are
incident with $v$.
For any $X,Y\subseteq E$, $X$ and $Y$
are called
{\it switching-equivalent,}
denoted by $X\ssim_G Y$,
if $X=Y\oplus \nabla_G(V_0)$ holds for a set $V_0$
of vertices in $G$,
where $A\oplus B$ is
the symmetric difference of two sets $A$ and $B$,
i.e., $A\oplus B=(A-B)\cup (B-A)$.
Let $X\nssim_G Y$ denote the case when
edge sets $X$ and $Y$ are not switching-equivalent in $G$.

Lukot'ka and Rollov\'a \cite{15}
proved that the property ``being feasible"
is invariant to switching-equivalent
edge sets.

\begin{thm}[\cite{15}]\relabel{SEP}
Let $G$ be a matching-covered graph
and $X$ and $Y$ be edge subsets of $G$.
If $X\ssim_G Y$,
then $X$ is feasible \iff $Y$ is  feasible.
\end{thm}

For a matching-covered graph $G=(V,E)$,
let $\F(G)$ be the set of feasible edge sets of $G$
and let $\nF(G)$ be the set of non-feasible edge sets of $G$.
Thus $\F(G)\cup \nF(G)$ is the power set of $E$.
Clearly $\{\emptyset, E\}\subseteq \nF(G)$.
Theorem~\ref{SEP} implies that
$\{X\subseteq E: X\ssim_G \emptyset\} \subseteq \nF(G)$
and $\{X\subseteq E: X\ssim_G  E\}\subseteq \nF(G)$.
For bipartite and regular graphs,
Lukot'ka and Rollov\'a~\cite{15} got the following
conclusion, described by notations in this article.

\begin{thm}[\cite{15}]
\relabel{luk}
If $G$ is a bipartite and regular graph, then
$\nF(G) = \{X\subseteq E: X\ssim_G \emptyset\}$.
\end{thm}

Note that any bipartite and regular graph is matching-covered,
because any bipartite graph is a class 1 graph (see \cite{Konig})
and any regular graph of class 1 is matching-covered.

In this article, we will extend Theorem~\ref{luk}
as stated below.

\begin{thm}\relabel{main1}
Let $G=(V,E)$ be a matching-covered graph.
Then the following statements are equivalent:
\begin{enumerate}
\item $G$ is bipartite;
\item
$\nF(G) = \{X\subseteq E: X\ssim_G \emptyset\}$;
\item $\nF(G) = \{X\subseteq E: X\ssim_G E\}$.
\end{enumerate}
\end{thm}

For any matching-covered graph $G=(V,E)$,
$\{X\subseteq E: X\ssim_G \emptyset
\mbox{ or } X\ssim_G  E\}$
is a subset of $\nF(G)$.
Let $\nF^*(G)=\nF(G)-\{X\subseteq E: X\ssim_G \emptyset
\mbox{ or } X\ssim_G  E\}$.
Then $\nF^*(G)=\emptyset$ holds \iff
$X\ssim_G \emptyset$ or
$X\ssim_G  E$ holds for each $X\in \nF(G)$.
It is natural to ask when $\nF^*(G)=\emptyset$ holds.
By Theorem~\ref{main1}, it holds if $G$ is bipartite
and matching-covered.
But there exist non-bipartite matching-covered graphs
with this property.
For example, $K_4$ is such a graph.

For a subgraph $G'$ of $G$, a {\it single ear} of $G'$
is a path $P$ of $G$ with an odd length
such that both ends of $P$ are in $G'$
but its internal vertices are distinct from
vertices in $G'$.
A {\it double ear} of $G'$ is a pair of
vertex disjoint single ears of $G'$.
An ear of $G'$ means a single ear or a double ear of $G'$.
An {\it ear decomposition}
of a matching-covered graph $G$ is a sequence
$$
G_0\subset G_1\subset\cdots\subset G_r=G
$$
of matching-covered subgraphs of $G$,
where (i) $G_0=K_2,$ and (ii)
for each $i$ with $1\leq i\leq r$,
$G_{i}$ is the union of $G_{i-1}$ and
an ear (single or double) of $G_{i-1}$.
For $i=1,2,\cdots,r$,
let $\epsilon(G_{i-1},G_i)\in \{1,2\}$
such that $\epsilon(G_{i-1},G_i)=1$
\iff $G_i$ is the union $G_{i-1}$ and a single ear.
A very important result on the study of matching-covered
graphs is the existence of an ear decomposition
for each matching-covered graph due to
Lov\'asz and Plummer \cite{19}.

Our second aim in this article is to
establish the following conclusions
on matching-covered graphs $G$
with $\nF^*(G)=\emptyset$,
based on ear decompositions of matching-covered graphs.

\begin{thm}\relabel{main2}
Let $G=(V,E)$ be a matching-covered graph
with an ear decomposition
$G_0\subset G_1\subset \cdots \subset G_r=G$,
where $r\ge 1$.
\begin{enumerate}
\item\relabel{main2-n1} If $\nF^*(G_{r-1})=\emptyset$ and
$\epsilon(G_{r-1},G_r)=1$,
then $\nF^*(G)=\emptyset$;
\item\relabel{main2-n2} if $\sum_{1\le i\le r}\epsilon(G_{i-1},G_i)\le r+1$,
then $\nF^*(G)=\emptyset$ holds;
\item\relabel{main2-n3} if
$\sum_{1\le i\le r}\epsilon(G_{i-1},G_i)\ge r+2$
and $\epsilon(G_{r-1},G_r)=2$,
then $\nF^*(G)\ne \emptyset$;
\item\relabel{main2-n4} if $\sum_{1\le i\le r}\epsilon(G_{i-1},G_i)\ge r+2$
and $\epsilon(G_{r-1},G_r)=1$, then
$\nF^*(G)= \emptyset$
\iff
$X\cap E(G_{r-1}-\{u,v\})$ is feasible in
the subgraph $G_{r-1}-\{u,v\}$
for each $X\in \nF^*(G_{r-1})$,
where $E(H)$ is the edge set of a graph $H$
and $u,v$ are the two ends  of the single ear $P_r$
added to $G_{r-1}$ for obtaining $G_r$.
\end{enumerate}
\end{thm}

Note that the graph $G_{r-1}-\{u,v\}$
in Theorem~\ref{main2} ~\ref{main2-n4}
is the graph obtained from $G_{r-1}$ by
deleting $u$ and $v$
and may not be matching-covered
although it contains perfect matchings.
By  definition, $X'=X\cap E(G_{r-1}-\{u,v\})$ is feasible in
$G_{r-1}-\{u,v\}$ if there exist two perfect matchings
$N_1$ and $N_2$ in $G_{r-1}-\{u,v\}$ such that
$|N_1\cap X'|\not \equiv  |N_2\cap X'|\modtwo$ holds.

Lukot'ka and Rollov\'a \cite{15}
noticed that $\nF^*(P)\ne \emptyset$ holds
for the Petersen graph $P$,
which is a class 2 graph,
and asked the following problem on regular graphs of class 1,
described by notations in this article.

\begin{prob}\relabel{prob1}
Does $\nF^*(G)=\emptyset$ hold for
each regular graph $G$ of class 1?
\end{prob}

A negative answer to this problem was provided by
He, et al \cite{wei} who showed that
for any $k\ge 3$, there exist infinitely many $k$-regular
graphs $G$ of class 1
with an arbitrary large equivalent edge set belonging to $\nF^*(G)$,
where a non-empty edge set $S$ of $G$
is called an {\it equivalent set}
if  $S\cap M=\emptyset$ or  $S\cap M=S$ holds
for all perfect matchings $M$ of $G$.
The graphs constructed in \cite{wei}
giving a negative answer to Problem~\ref{prob1}
are not 3-connected and the following problem was further asked
in \cite{wei}.

\begin{prob}\relabel{prob2}
Does Problem~\ref{prob1} hold
for $3$-connected and
$r$-regular graph $G$ with $r\ge 3$?
\end{prob}

In Section~\ref{sec5}, we will provide negative answers
to both Problems~\ref{prob1} and~\ref{prob2} by
two constructions of
$r$-regular graphs $G$ of class 1 with $\nF^*(G)\ne \emptyset$.

\begin{thm}\relabel{main3}
For any integer $r\ge 3$,
there are infinitely many $r$-connected and
$r$-regular graphs $G$ of class 1
with $\nF^*(G)\ne \emptyset$.
\end{thm}

\section{Preliminary results on
$X\subseteq E$ with $X\ssim_G \emptyset$
or $X\ssim_G E$
\relabel{sec2}}

Let $G=(V,E)$ be any connected graph
which may be not matching-covered.
By definition,
for any subset $U\subseteq V$,
$\nabla_G(U)$ is  the set
$\{e\in E: e$
joins a vertex in in  $U$
and a vertex in $V-U\}$.
With the notation $\nabla_G(U)$,
an edge set $X$ of $G$ with the property that
$X\ssim_G \emptyset$
or $X\ssim_G E$ has the following
characterization due to He, et al \cite{wei}.

\begin{pro}[\cite{wei}]\relabel{pro2-1}
Let $G=(V,E)$ be a connected graph and $X\subseteq E$.
Then
\begin{enumerate}
\item $X\ssim_G \emptyset$  iff
$X=\nabla _G(U)$ for some $U\subseteq V$;
\item $X\ssim_G E$  \iff
$E(G)-X=\nabla _G(U)$ for some $U\subseteq V$.
\end{enumerate}
\end{pro}

Proposition~\ref{pro2-1} implies the following corollary
immediately.
For any graph $G$ and any set $V_0$ of vertices in $G$,
let $G[V_0]$ denote the subgraph of $G$ induced by $V_0$.

\begin{cor}\relabel{cor2-1}
Let $G=(V,E)$ be a connected graph and $X\subseteq E$.
For any $V_0\subseteq V$,
\begin{enumerate}
\item
if $X\ssim_G \emptyset$, then
$X\cap E(G[V_0])\ssim_{G[V_0]} \emptyset$;
\item
if $X\ssim_G E$, then
$X\cap E(G[V_0])\ssim_{G[V_0]} E(G[V_0])$.
\end{enumerate}
\end{cor}

Obviously,
$X\ssim_G Y$ implies that $Y\ssim_G X$.
The transitive property of the relation ``$\ssim_G$" also holds.

\begin{lem}\relabel{le2-0}
Let $G=(V,E)$ be a connected graph with $X,Y,Z\subseteq E$.
If $X\ssim_G Y$ and $Y\ssim_G Z$, then
$X\ssim_G Z$ holds.
\end{lem}

\proof Assume that
$X\ssim_G Y$ and $Y\ssim_G Z$.
Then $X=Y\oplus \nabla_G(V_1)$ and
$Y=Z\oplus \nabla_G(V_2)$ hold for some
$V_1,V_2\subseteq V$,
implying that $X=Z\oplus \nabla_G(V_1\oplus V_2)$.
Thus  $X\ssim_G Z$ holds.
\endproof

\vspace{0.3 cm}

Assume that $G'$ is any connected graph with
two distinct vertices $v_1$ and $v_2$
and $P$ is any path with ends $u_1$ and $u_2$
such that $G'$ and $P$ are vertex-disjoint.
Let $\Union{(v_1,v_2)}{G',P}$
(or simply $\Union{}{G',P}$)
denote the graph obtained
from $G'$ and $P$ by identifying $u_i$ and $v_i$
for $i=1,2$.
For an ear decomposition $G_0\subset G_1\subset \cdots G_r=G$
of a matching-covered graph $G$,
if $G_i$ is the union of $G_{i-1}$ and a single ear $P_i$,
then $G_i=\Union{}{G_{i-1},P_i}$.
But, in this section, the results do not depend on
the condition that $G'$ is matching-covered.

\begin{lem}\relabel{le2-1}
Let $G=\Union{}{G',P}$.
For any edge set $X=X_0\cup X'$ of $G$,
where $X_0\subseteq E(P)$ and $X'\subseteq E(G')$,
\begin{enumerate}
\item\relabel{le2-1-n1}
if $|E(P)|\equiv 1\modtwo$ and $X\in \nF(G)$,
then $X'\in \nF(G')$;
\item\relabel{le2-1-n2} if $|X_0|\equiv 0\modtwo$, then $X\ssim_G X'$;
\item\relabel{le2-1-n3} if $|X_0|\equiv 1\modtwo$,
then $X\ssim_G X'\cup \{e\}$
for any $e\in E(P)$;
\item\relabel{le2-1-n4} if $X'\ssim_{G'} Y$, then $X\ssim_G Y\cup Y_0$
for some $Y_0\subseteq E(P)$;
\item\relabel{le2-1-n5} if $X'\ssim_{G'} \emptyset$,
then either $X\ssim_{G} \emptyset$ or
$X\ssim_{G} \{e\}$ for any $e\in E(P)$;
\item\relabel{le2-1-n6} if $X'\ssim_{G'} E(G')$,
then either $X\ssim_{G} E(G)$ or
$X\ssim_{G} E(G)-\{e\}$ for any $e\in E(P)$;
\item\relabel{le2-1-n7}
if $X\ssim_{G} \emptyset$, then $X'\ssim_{G'} \emptyset$;
if $X\ssim_{G} E(G)$, then $X'\ssim_{G'} E(G')$.
\end{enumerate}
\end{lem}

\proof
(i). Assume that the edges in $P$
are $e_1, e_2,\cdots, e_{2k-1}$ in
the order of the path $P$ such that $e_i$ and $e_{i+1}$
have a common end for all $i=1,2,\cdots,2k-2$.
Suppose that $X'\in \F(G')$.
Then $G'$ has two perfect matchings $M_1$ and $M_2$
such that
$|X'\cap M_1|-|X'\cap M_2|\equiv 1\modtwo$.
For $i=1,2$, the set $N_i$ defined below is a
perfect matching of $G$:
$$
N_i=M_i\cup
\{e_{2j}: j=1,2,\cdots, k-1\}.
$$
Observe that
$$
|X\cap N_j|=|X'\cap M_j|+|X_0\cap
\{e_{2j}: j=1,2,\cdots, k-1\}|,\qquad \forall j=1,2.
$$
Thus $|X\cap N_1|-|X\cap N_2|
=|X'\cap M_1|-|X'\cap M_2|
\equiv 1\modtwo$,
implying that
$X$ is feasible in $G$, a contradiction.

Thus \ref{le2-1-n1} holds.

\ref{le2-1-n2} and \ref{le2-1-n3} will be proved
by applying the following claim.

\noindent {\bf Claim 1}:
If $|X_0|\ge 2$, then $X=X_0\cup X'\ssim_G X_0'\cup X'$ holds
for some $X_0'\subset X_0$ with $|X'_0|=|X_0|-2$.

Assume that $|X_0|\ge 2$.
Then there exists subpath $P_0$ of $P$
such that $X\cap E(P_0)=\emptyset$
and $\nabla_G(V(P_0))\subseteq X$,
implying that $X\ssim_G X\oplus \nabla_G(V(P_0))
=X'_0\cup X'$,
where $X'_0=X_0\oplus \nabla_G(V(P_0))\subset X_0$
and $|X'_0|=|X_0|-2$.
Thus the claim holds.

(ii).
Assume that $|X_0|>0$ and $|X_0|\equiv 0\modtwo$.
\ref{le2-1-n2}
follows by applying Claim 1
repeatedly.

(iii). Applying Claim 1 repeatedly,
$X\ssim_G \{e\}\cup X'$ holds for some $e\in E(P)$.
Now let $e'$ be any edge in $P$ different from $e$.
There exists a subpath $P'$ of $P$
such that $\nabla_G(V(P'))=\{e,e'\}$.
Thus $\{e\}\cup X'\ssim_G
(\{e\}\cup X')\oplus \nabla_G(V(P'))=
\{e'\}\cup X'$
and the result holds.

(iv). It is trivial when $X'=Y$.
Now assume that $X'\ne Y$. Then
$Y=X'\oplus \nabla_{G'}(V_0)$ for some non-empty
set $V_0\subset V(G')$.

As $G=\Union{}{G',P}$,
there are three cases on the structure of $G$, i.e., $|\{v_1,v_2\}\cap V_0|\in \{0,1,2\}$,
where $v_1,v_2$ are the two vertices in $G'$
at which the ends of $P$ are identified with.
But $|\{v_1,v_2\}\cap V_0|=2$ implies that
$|\{v_1,v_2\}\cap (V(G')-V_0)|=0$.
Thus, we need only to consider the two cases:
$|\{v_1,v_2\}\cap V_0|=0$ or
$|\{v_1,v_2\}\cap V_0|=1$,
as shown in Figure~\ref{f1}.
\begin{figure}[ht]
\centering
\input{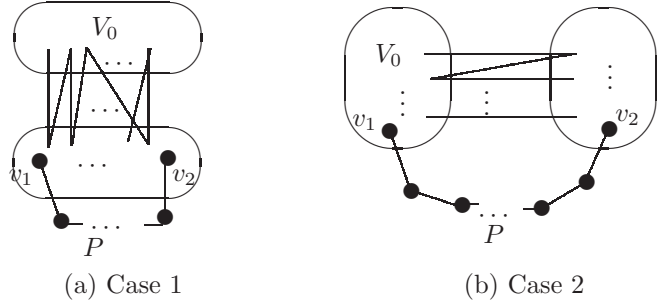}

(a) Case 1 \hspace{3.5 cm} (b) Case 2

\caption{Two cases for the two ends of $P$}
\relabel{f1}
\end{figure}

In both cases, $Y=X'\oplus \nabla_{G'}(V_0)$ implies that
$X\oplus \nabla_{G}(V_0) =Y\cup Y_0$ holds for some
$Y_0\subseteq E(P)$.

Thus \ref{le2-1-n4} holds.

(v). As $X'\ssim_{G'}\emptyset$,
the result of \ref{le2-1-n4} implies that
$X\ssim_{G} Y_0$ where $Y_0\subseteq E(P)$.
The results of \ref{le2-1-n2} and \ref{le2-1-n3} imply that
either $Y_0\ssim_G \emptyset$ or
$Y_0\ssim_G \{e\}$ for any $e\in E(P)$.

Thus \ref{le2-1-n5} holds.

(vi).
As $X'\ssim_{G'} E(G')$,
the result of \ref{le2-1-n4} implies that
$X\ssim_{G} E(G')\cup Y_0$ where $Y_0\subseteq E(P)$.
The results of \ref{le2-1-n2} and \ref{le2-1-n3} imply that
$(E(G')\cup Y_0)\ssim_G E(G)$
when $|Y_0|\equiv |E(P)|\modtwo$,
and  $(E(G')\cup Y_0)\ssim_G E(G)-\{e\}$
when $|Y_0|\not\equiv |E(P)|\modtwo$.

Thus \ref{le2-1-n6} holds.

(vii).
Suppose that $X\ssim_G \emptyset$ holds.
By Proposition~\ref{pro2-1},
$X=\nabla_G(U)$ holds for some $U\subseteq V(G)$.
Then $X'=\nabla_{G'}(U-U_0)$,
implying that $X'\ssim_{G'} \emptyset$,
where $U_0$ is the set of internal vertices of $P$.

Now suppose that $X\ssim_G E(G)$ holds.
By Proposition~\ref{pro2-1},
$E(G)-X=\nabla_G(U)$ holds for some $U\subseteq V(G)$.
Then $E(G')-X'=\nabla_{G'}(U-U_0)$,
where $U_0$ is defined above,
implying that $X'\ssim_{G'} E(G')$.

Thus  \ref{le2-1-n7} holds.
\endproof

\vspace{0.3 cm}

For distinct vertices $v_1,v_2,v_3,v_4$
in a graph $G'$ and any two vertex-disjoint
paths $P_1,P_2$ with $V(P_i)\cap V(G')=\emptyset$
for $i=1,2$, let
$\Union{(v_1,v_2,v_3,v_4)}{G',P_1,P_2}$
(or simply $\Union{}{G',P_1,P_2}$)
be the graph $\Union{(v_3,v_4)}{G'',P_2}$,
where $G''=\Union{(v_1,v_2)}{G',P_1}$.

\begin{lem}\relabel{le2-2}
Let $G=\Union{}{G',P_1,P_2}$.
For any edge set $X=X_0\cup X'$ of $G$,
where $X'\subseteq E(G')$
and $X_0\subseteq E(P_1)\cup E(P_2)$,
if $X'\ssim_{G'} \emptyset$,
then $X\ssim_G \emptyset$,
or $X\ssim_G \{e\}$ for some $e\in E(P_1)\cup E(P_2)$,
or $X\ssim_G \{e_1,e_2\}$ where $e_i\in E(P_i)$
for $i=1,2$.
\end{lem}

\proof Let $G''=\Union{}{G',P_1}$.
As $X'\ssim_{G'} \emptyset$,
Lemma~\ref{le2-1}~\ref{le2-1-n5}
implies that $X-E(P_2)\ssim_{G''} \emptyset$
or  $X-E(P_2)\ssim_{G''} \{e_1\}$ for any $e_1\in E(P_1)$.

Note that $G=\Union{}{G'',P_2}$.
If $X-E(P_2)\ssim_{G''} \emptyset$,
then Lemma~\ref{le2-1}~\ref{le2-1-n5}
implies that either $X\ssim_{G} \emptyset$
or  $X\ssim_{G} \{e_2\}$ holds for any $e_2\in E(P_1)$.

If $X-E(P_2)\ssim_{G''} \{e_1\}$ for any $e_1\in E(P_1)$,
then Lemma~\ref{le2-1}~\ref{le2-1-n4}
implies that  $X\ssim_{G} \{e_1\}\cup Y_0$
for some $Y_0\subseteq E(P_2)$.
Lemma~\ref{le2-1}~\ref{le2-1-n2} and~\ref{le2-1-n3}
further imply that either $X\ssim_{G} \{e_1\}$  or
$X\ssim_{G} \{e_1\}\cup \{e_2\}$ holds
for any $e_2\in E(P_2)$.

Thus the result holds.
\endproof

\section{Proof of Theorem~\ref{main1}\relabel{sec3}}

For any matching-covered graph $G$, the following
basic properties follow directly from
the definitions of $\F(G)$ and $\nF(G)$.

\begin{lem}\relabel{le3-1}
Let $G$ be a matching-covered graph with $|E(G)|\ge 2$
and $X\subseteq E(G)$.
If  either $|X|=1$ or $|X|=|E|-1$,
then $X\in \F(G)$.
\end{lem}

By applying Lemma~\ref{le2-1}, we can prove that
for any matching-covered graphs $G'$ and  $G=\Union{}{G',P}$,
$\nF^*(G')=\emptyset$ implies that
$\nF^*(G)=\emptyset$.

\begin{lem}\relabel{le3-2}
Let $G'$ and
$G=\Union{}{G',P}$ be matching-covered graphs,
where $P$ is a single ear of $G'$.
For any $X\in \nF(G)$,
\begin{enumerate}
\item $X\cap E(G')\ssim_{G'} \emptyset$
\iff $X\ssim_G \emptyset$;
\item $X\cap E(G')\ssim_{G'} E(G')$ \iff  $X\ssim_G E(G)$;
\item $X\cap E(G')\in \nF^*(G')$
\iff $X\in \nF^*(G)$.
\end{enumerate}
\end{lem}

\proof (i). ($\Leftarrow$) It follows directly from
Lemma~\ref{le2-1}~\ref{le2-1-n7}.

($\Rightarrow$)
As $X\cap E(G')\ssim_{G_{r-1}} \emptyset$,
Lemma~\ref{le2-1}~\ref{le2-1-n5} implies that
$X\ssim_G \emptyset$ or $X\ssim_G \{e\}$ for any $e\in E(P_r)$.

Suppose that $X\ssim_G \{e\}$ for some $e\in E(P_r)$.
As $X\in \nF(G)$, Theorem~\ref{SEP}
implies that $\{e\}\in \nF(G)$.
But, as $|E(G)|\ge 2$,
Lemma~\ref{le3-1}
implies that $\{e\}\in \F(G)$,
a contradiction.
Thus $X\ssim_G \emptyset$.

(ii). ($\Leftarrow$) It follows directly from
Lemma~\ref{le2-1}~\ref{le2-1-n7}.

($\Rightarrow$).
As $X\cap E(G')\ssim_{G_{r-1}} E(G')$,
Lemma~\ref{le2-1}~\ref{le2-1-n6} implies that
$X\ssim_G E(G)$ or $X\ssim_G E(G)-\{e\}$ for any $e\in E(P_r)$.

Suppose that $X\ssim_G E(G)-\{e\}$ for some $e\in E(P_r)$.
As $X\in \nF(G)$, Theorem~\ref{SEP}
implies that $E(G)-\{e\}\in \nF(G)$ holds.
As $|E(G)|\ge 2$,
Lemma~\ref{le3-1}
implies that $E(G)-\{e\}\in \F(G)$,
a contradiction.
Thus $X\ssim_G E(G)$.

(iii). By Lemma~\ref{le2-1}~\ref{le2-1-n1},
$X\in \nF(G)$ implies that $X\cap E(G')\in \nF(G')$.
Then the result follows from (i) and (ii) directly.
\endproof

\vspace{0.3 cm}

An ear decomposition
$G_0\subset G_1\subset \cdots \subset G_r$
of a matching-covered graph $G$ is called
a {\it single-ear decomposition}
if $\epsilon(G_{i-1},G_i)=1$ holds for all
$i=1,\cdots,r$.
A matching-covered graph may have no
single-ear decompositions.
For example, the complete graph $K_4$ does not have.
However, every matching-covered bipartite graph
has a single-ear decomposition.

\begin{thm}
\relabel{Ear-c}
Let $G$ be a matching-covered  graph $G$.
\begin{enumerate}
\item\ \cite{19} $G$ has an ear decomposition;
\item\ \cite{Ear,L-Ear,8}
$G$ is bipartite \iff
$G$ has a single-ear decomposition.
\end{enumerate}
\end{thm}

Now we are going to prove Theorem~\ref{main1}.

\vspace{0.4 cm}

{\bf Proof of Theorem~\ref{main1}}:
By Proposition~\ref{pro2-1} (i),
each edge set $X$ with $X\ssim_G \emptyset$
induces a bipartite subgraph in $G$,
implying that $E\nssim_G \emptyset$ holds whenever
$G$ is not bipartite.
Hence, Theorem~\ref{main1} (ii) implies
Theorem~\ref{main1} (i).

For any bipartite graph $G=(V,E)$,
$E\ssim_G \emptyset$ holds.
Thus, Lemma~\ref{le2-0} implies that
$\{X\subseteq E: X\ssim_G \emptyset\}$
and $\{X\subseteq E: X\ssim_G E\}$
are the same set.
Thus,
(ii) and (iii) in Theorem~\ref{main1} are equivalent.

So, to prove Theorem~\ref{main1},
it suffices to show that
Theorem~\ref{main1} (i) implies Theorem~\ref{main1} (ii).

Assume that $G$ is bipartite and matching-covered.
By Theorem \ref{Ear-c} (ii),
$G$ has a single-ear decomposition
$G_0\subset G_1\subset\cdots\subset G_r=G$,
where $G_0\cong K_2$.
Thus, for $i=1,2,3,\cdots,r$,
$G_i=\Union {}{G_{i-1},P_i}$ holds for some single ear $P_i$
of $G_{i-1}$.

If $r=0$, i.e., $G\cong K_2$
and (i) implies (ii) obviously.
Now assume that $r\ge 1$ and the result holds for $G_{r-1}$.
For any $X\in \nF(G)$,
Lemma~\ref{le2-1}~\ref{le2-1-n1} implies that
$X\cap E(G_{r-1})\in \nF(G_{r-1})$.
By the assumption,  the result holds for $G_{r-1}$.
Thus  $X\cap E(G_{r-1})\ssim_{G_{r-1}} \emptyset$ holds.
Then, Lemma~\ref{le3-2} (i) implies that
$X\ssim_G \emptyset$.

Hence Theorem~\ref{main1} is proven.
\endproof

\section{Proof of Theorem~\ref{main2}\relabel{sec4}}

The following two  lemmas will be applied for
proving Theorem~\ref{main2}.

\begin{lem}\relabel{le4-1}
Let $G'$ and
$G=\Union{}{G',P_1,P_2}$ be matching-covered graphs,
where $P_1$ and $P_2$ form a double ear of $G'$.
Assume that $\Union{}{G',P_i}$ is not matching-covered
for $i=1,2$.
Then $\nF^*(G)=\emptyset$ \iff $G'$ is bipartite.
\end{lem}

\proof ($\Rightarrow $)
Suppose that $G'$ is not bipartite.

Let $X_0=E(P_1)\cup E(P_2)$, where
$E(P_i)=\{e_{i,1},e_{i,2},\cdots, e_{i,2k_i-1}\}$
for $i=1,2$ and $e_{i,j}$ and $e_{i,j+1}$ have a common end
for all $j=1,2,\cdots,2k_i-2$.
As  $\Union{}{G',P_i}$ is not matching-covered
for both $i=1,2$,
for each perfect matching $M$ of $G$,
one of the following holds:
$$
M\cap X_0=
\bigcup_{1\le i\le 2}
\{e_{i,2t-1} : t=1,2,\cdots,k_i\}
$$
or
$$
M\cap X_0=
\bigcup_{1\le i\le 2}
\{e_{i,2t} : t=1,2,\cdots,k_i-1\}.
$$
Thus $|M\cap X_0|\equiv 0\modtwo$
holds for all perfect matchings $M$ of $G$,
implying that $X_0\in \nF(G)$.

As $G'$ is not bipartite,
$G-X_0$ is not bipartite.
Thus
Proposition~\ref{pro2-1} (ii) implies that
$X_0\nssim_{G} E(G)$.

As $|E(P_1)|\equiv |E(P_2)|\equiv 1\modtwo$,
Lemma~\ref{le2-1}~\ref{le2-1-n3} implies that
$X_0\ssim_G \{e_1,e_2\}$,
where $e_i$ is an edge on $P_i$ for $i=1,2$.
Clearly $G-\{e_1,e_2\}$ is connected.
Then Proposition~\ref{pro2-1} (ii) implies that
$\{e_1,e_2\}\nssim_{G} \emptyset$.
Thus $X_0\nssim_G \emptyset$.

Hence $X_0\in \nF^*(G)$ and the necessity holds.

($\Leftarrow $)
Assume that $G'$ is bipartite
and $X\in \nF(G)$.

As $|E(P_1)|\equiv |E(P_2)|\equiv 1\modtwo$,
Lemma~\ref{le2-1}~\ref{le2-1-n1} implies that
$X\cap E(G')\in \nF(G')$.
As $G'$ is bipartite,
Theorem~\ref{main1} implies that
$X\cap E(G')\ssim_{G'} \emptyset$.
By Lemma~\ref{le2-2},
$X\ssim_G \emptyset$ holds or
$X\ssim_G \{e\}$ holds for some $e\in E(P_1)\cup E(P_2)$
or
$X\ssim_G \{e_1,e_2\}$ holds for some $e_1\in E(P_1)$
and $e_2\in E(P_2)$.

If $X\ssim_G \{e\}$ for some $e\in E(P_1)\cup E(P_2)$,
then Theorem~\ref{SEP} implies that $\{e\}\in \nF(G)$.
But Lemma~\ref{le3-1} implies that $\{e\}\in \F(G)$,
a contradiction.

Now consider the case that
$X\ssim_G \{e_1,e_2\}$ for $e_i\in E(P_i)$.
Lemma~\ref{le2-1}~\ref{le2-1-n3} implies that
$\{e_1,e_2\}\ssim_G (E(P_1)\cup E(P_2))$.
Thus $X\ssim_G (E(P_1)\cup E(P_2))$.
Since $G'$ is bipartite and matching-covered, $G'$ has a bipartition $(U_1,U_2)$ with $|U_1|=|U_2|$.
Since $G=\Union{}{G',P_1,P_2}$ is not bipartite,
both ends of some $P_i$
are within $U_j$ for some $j$.
Assume that both ends of some $P_1$ are within $U_1$.
As $|U_1|=|U_2|$ and $G$ is matching-covered,
both ends of some $P_2$ must be in $U_2$.
Thus  $(E(P_1)\cup E(P_2))\oplus
\nabla_G(U_1\cup V(P_1))
=E(G)$,
implying that $X\ssim_G (E(P_1)\cup E(P_2)) \ssim_G E(G)$.

Hence the sufficiency holds.
\endproof

\begin{lem}\relabel{le4-2}
Let $G'$ and $G=\Union{}{G',P}$ be
matching-covered graphs,
where $P$ is a single ear of $G'$.
For any $X\in \nF(G')$,
both $X \in \F(G)$ and $X\cup E(P)\in \F(G)$ hold
\iff
$X\cap E(G^o)\in \F(G^o)$ holds, where $G^o=G'-\{u,v\}$
and $u,v$ are the two ends of $P$ in $G'$.
\end{lem}

\proof As $P$ is a single ear of $G'$, $|E(P)|$ is odd.
Let $e_1,e_2,\cdots,e_{2k-1}$ be the edges in $P$,
where $e_i$ and $e_{i+1}$ have a common end
for all $i=1,2,\cdots,2k-2$.

The set of perfect matchings of $G$ can be partitioned
into two sets $\M_0$ and $\M_1$, where
$\M_0$ is the set of perfect matchings $M$
in $G$ with $e_1\notin M$
and $\M_1$ is the set of perfect matchings $M$
in $G$ with $e_1\in M$.
Then, for each $M\in \M_0$,
$$
M\cap E(P)=\{e_{2r}: r=1,2,\cdots,k-1\}
$$
and
for each $M\in \M_1$,
$$
M\cap E(P)=\{e_{2r-1}: r=1,2,\cdots,k\}.
$$

Observe that $\M'=\{M\cap E(G'): M\in \M_0\}$
is the set of perfect matchings in $G'$.
Assume that $X\in \nF(G')$ and
$|M'\cap X|\equiv a\modtwo$ holds for all $M'\in \M'$,
where $a$ is a fixed number in $\{0,1\}$.
Thus $|M\cap X|\equiv a+k-1 \modtwo$ holds
for all $M\in \M_0$.

$(\Rightarrow )$ Assume that both $X$ and
$X\cup E(P)$ are feasible in $G$.
Since $X$ is feasible in $G$
and $|M\cap X|\equiv a+k-1 \modtwo$ holds
for all $M\in \M_0$,
$|M_1\cap X|\equiv a+k \modtwo$ holds
for some $M_1\in \M_1$.

{\bf Claim 1}:
$|M_2\cap X|\equiv a+k-1 \modtwo$ holds
for some $M_2\in \M_1$.

Suppose that Claim 1 fails.
Then $|M\cap X|\equiv a+k \modtwo$ holds
for all $M\in \M_1$,
implying that
$|M\cap (X\cup E(P))|\equiv a \modtwo$ holds
for all $M\in \M_1$.
But, for each $M\in M_0$,
$|M\cap (X\cup E(P))|\equiv |M\cap X|+k-1\equiv a \modtwo$ holds.
Thus $|M\cap (X\cup E(P))| \equiv a \modtwo$ holds
for all $M\in  \M_0\cup \M_1$,
implying that $X\cup E(P)$ is non-feasible in $G$, a contradiction.

Thus Claim 1 holds.

Now there are two perfect matchings $M_1,M_2\in \M_1$
such that $|M_i\cap X|\equiv a+k+i-1 \modtwo$ holds
for $i=1,2$,
implying that $|M_1\cap X|\not \equiv |M_2\cap X|\modtwo$.

Let $X_0=X\cap E(G^o)$.
Observe that both $M_1-E(P)$ and $M_2-E(P)$
are perfect matchings in $G^o$ and
$
|(M_i-E(P))\cap X_0|=|M_i\cap X|
$
holds for $i=1,2$.
As $|M_1\cap X|\not \equiv |M_2\cap X|\modtwo$,
$|(M_1-E(P))\cap X_0|\not\equiv
|(M_2-E(P))\cap X_0|\modtwo$ holds,
implying that $X_0$ is feasible in $G^o$.

$(\Leftarrow )$
Assume that $X_0=X\cap E(G^o)$ is feasible in $G^o$.
Then there are two perfect matchings
$N_1$ and $N_2$ in $G^o$ such that
$|X_0\cap N_i|\equiv i \modtwo $ for $i=1,2$.

Clearly, $Q_i=N_i\cup \{e_{2r-1}: r=1,2,\cdots,k\}\in \M_1$
 for $i=1,2$.
Observe that
$$
|(X\cup E(P))\cap Q_i|=|X_0\cap N_i|+k
\equiv k+i \modtwo,\quad \forall i=1,2,
$$
 implying that $X\cup E(P)$ is feasible in $G$.
Also observe that
$$
|X\cap Q_i|=|X_0\cap N_i|
\equiv i \modtwo, \quad \forall i=1,2,
$$
implying that $X$ is feasible in $G$.
\endproof

\vspace{0.4 cm}

We are now ready to prove Theorem~\ref{main2}.

\vspace{0.2 cm}

{\bf Proof of Theorem~\ref{main2}}:
(i). It follows directly from Lemma~\ref{le3-2} (iii).

(ii). If $\sum_{1\le i\le r}\epsilon(G_{i-1},G_i)=r$,
then $G_0\subset G_1\subset \cdots \subset G_r$ is a single
ear decomposition of $G$. Thus Theorem~\ref{main1}
implies that $\nF^*(G)=\emptyset$.

Now assume that $\sum_{1\le i\le r}\epsilon(G_{i-1},G_i)=r+1$,
implying that
$\epsilon(G_{i-1},G_i)=2$ holds for exactly one $i$
with $1\le i\le r$.
We first consider the case that $\epsilon(G_{r-1},G_r)=2$.
In this case,
$\sum_{1\le i\le r-1}\epsilon(G_{i-1},G_i)=r-1$,
implying that
$G_0\subset G_1\subset \cdots \subset G_{r-1}$ is a single
ear decomposition of $G_{r-1}$.
Theorem~\ref{Ear-c} implies that $G_{r-1}$ is bipartite.
Then Lemma~\ref{le4-1} implies that
$\nF^*(G_r)=\emptyset$ holds.

Now we consider the case
that $\epsilon(G_{k-1},G_k)=2$, where $1\le k<r$.
Then $\sum_{1\le i\le k}\epsilon(G_{i-1},G_i)=k+1$.
By the proven conclusion above,
$\nF^*(G_k)=\emptyset$ holds.
The result in (i) implies that
$\nF^*(G_k)=\emptyset$ holds
for all $i=k+1,\cdots,r$.

Hence (ii) holds.

(iii). As $\sum\limits_{1\le i\le r}\epsilon(G_{i-1},G_i)\ge r+2$
and $\epsilon(G_{r-1},G_r)=2$,
$\sum\limits_{1\le i\le r-1}\epsilon(G_{i-1},G_i)\ge r$ holds.
By the definition of ear decompositions,
$G_{r-1}$ is not bipartite.
Then Lemma~\ref{le4-1} implies that $\nF^*(G)\ne \emptyset$.
Hence (iii) holds.

(iv).
($\Rightarrow $)
Assume that $\nF^*(G)=\emptyset$.
Suppose that there exists
$X\in \nF^*(G_{r-1})$ with
$X\cap E(G^o)\in \nF(G^o)$,
where $G^o=G_{r-1}-\{u,v\}$.

As $X\cap E(G^o)\in \nF(G^o)$,
Lemma~\ref{le4-2} implies that
$X\in \nF(G)$ or $X\cup E(P)\in \nF(G)$ holds.
If $X\in \nF(G)$,
as $X\in \nF^*(G_{r-1})$,
then Lemma~\ref{le3-2} (iii) implies that
$X\in \nF^*(G)$.
If $X\cup E(P)\in \nF(G)$,
it can be proved similarly that
$X\cup E(P)\in \nF^*(G)$ holds.
Thus $\nF^*(G)\ne \emptyset$,
a contradiction.

($\Leftarrow $)  Assume that $\nF^*(G)\ne \emptyset$.
Then, there exists $Z\in \nF^*(G)$ and Lemma~\ref{le3-2} (iii) implies that
$X=Z\cap E(G_{r-1})\in \nF^*(G_{r-1})$.

By Lemma~\ref{le2-1}~\ref{le2-1-n2} and~\ref{le2-1-n3},
$Z\ssim_G X$ or $Z\ssim_G X\cup E(P_r)$ holds.
Then, $Z\in \nF(G)$ implies that
$X\in \nF(G)$ or $X\cup E(P)\in \nF(G)$.
Lemma~\ref{le4-2} implies that $X\cap E(G^o)\in \nF(G^o)$ holds,
where $G^o=G_{r-1}-\{u,v\}$, contradicting the given condition.

Thus the result holds.
\endproof

\section{
Regular graphs $G$ of class 1 with
$\nF^*(G)\ne \emptyset$
\relabel{sec5}
}

\subsection{Generalize the family of graphs
constructed in \cite{wei}
\relabel{sec5-1}}

In this subsection, we will
generalize the construction
in \cite{wei} which provides a negative answer
to Problem~\ref{prob1}.

For two vertex-disjoint
graphs $G_1=(V_1,E_1)$ and $G_2=(V_2,E_2)$
with $e_i=x_iy_i\in E_i$ for $i=1,2$,
let $G_1\#_{e_1,e_2} G_2$ denote the graph obtained from $G_1-e_1$ and $G_2-e_2$
by adding edges $f_1=x_1x_2$ and $f_2=y_1y_2$,
as shown in Figure~\ref{f2}.

\begin{figure}[ht]
\centering
\input{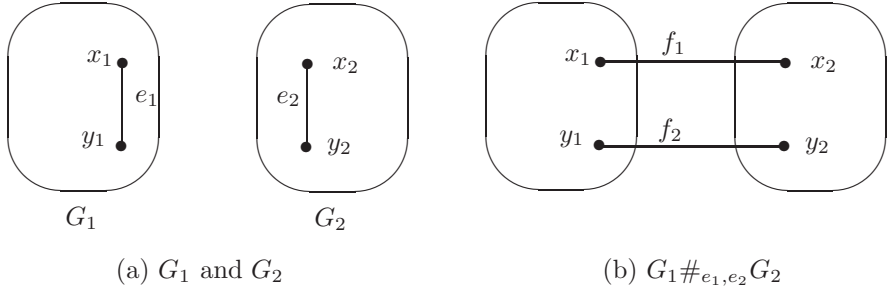}

(a) $G_1$ and $G_2$ \hspace{4 cm}
(b) $G_1\#_{e_1,e_2} G_2$

\caption{A graph constructed from $G_1$
and $G_2$}
\relabel{f2}
\end{figure}

\begin{lem}\relabel{le5-1-1}
For $i=1,2$, assume that
$G_i=(V_i,E_i)$ is a matching-covered graph
with $|E_i|\ge 2$ and
$S_i$ is an equivalent set of $G_i$ with $e_i\in S_i$,
where $e_i=x_iy_i$.
Let $G$ denote the graph $G_1\#_{e_1,e_2} G_2$
and let $S=(S_1-\{e_1\})\cup (S_2-\{e_2\})\cup \{f_1, f_2\}$.
Then
\begin{enumerate}
\item $G$ is
matching-covered;
\item
$S$
is an equivalent set in $G$;
\item when $G_1$ and $G_2$ are 2-connected,
$G$ is also 2-connected;

\item when $G_1$ and $G_2$ are $r$-regular graphs
of class 1,
$G$ is also a $r$-regular graph of class 1;
\item
for any $S'\subseteq S$,
when $G_i-e_i-(S'\cap E_i)$ is not bipartite
for some $i\in \{1,2\}$,
$S'\nssim_G E(G)$ holds;
\item
for any $S'\subseteq S$, when $S'\cap E_j \nssim_{G_j-e_j} \emptyset$ for some $j\in \{1,2\}$, $S'\nssim_G \emptyset$ holds.
\end{enumerate}
\end{lem}

\proof For $i=1,2$ and $j=0,1$,
let $\M_{i,j}$ be the set of perfect matchings $M$ in $G_i$
with $|M\cap \{e_i\}|=j$.
Since $G_i$ is matching-covered and
$|E_i|\ge 2$ holds for $i=1,2$,
$\M_{i,j}\ne \emptyset$ for all $i=1,2$ and $j=0,1$.
Let $\M$ be the set of perfect matchings of $G$.

(i).
The following facts imply that $G$ is matching-covered:
\begin{enumerate}
\item[(a)] for any $M_i\in \M_{i,1}$, $i=1,2$,
$(M_1-\{e_1\}) \cup (M_2-\{e_2\})\cup \{x_1y_1,x_2y_2\}$
is a member in $\M$;
\item[(b)] if $N_i\in \M_{i,0}$ for $i=1,2$,
then $N_1\cup N_2\in \M$;
\item[(c)] for $i=1,2$ and any $e\in E_i$,
$e\in M_i$ holds for some $M_i\in \M_{i,1}\cup \M_{i,0}$.
\end{enumerate}

(ii). To show that $S$ is an equivalent set of $G$,
we need only to prove the two claims below:

\noindent {\bf Claim 1}: For $\{f_1,f_2\}$
is an equivalent set of $G$.

Suppose the claim fails.
Then there exists $M\in \M$ with
$|\{f_1,f_2\}\cap M|=1$.
Assume that $f_1\in M$ but $f_2\notin M$.
Then $M\cap E_1$ is a perfect matching of $G-x_1$,
implying that $|V(G_1)|\equiv 1\modtwo$,
contradicting the condition that $G_1$ is matching-covered.
Thus the claim holds.

\noindent {\bf Claim 2}:
both $\{f_1, e\}$
is an equivalent set of $G$ for any
$e\in (S_1-\{e_1\})\cup (S_2-\{e_2\})$.

We may assume that $e\in S_1-\{e_1\}$.
Suppose the claim fails.
Then there exists $M\in \M$ with
$|\{f_1,e\}\cap M|=1$.

If $e\in M$ but $f_1\notin M$, then Claim 1 implies that
$f_2\notin M$.
Thus $M_1=M\cap E_1\in \M_{1,0}$.
Clearly,  $e\in M_1$ but $e_1\notin M_1$.
Thus $\{e,e_1\}$ is not an equivalent set of $G_1$,
contradicting the assumption that
$S_1$ is an equivalent set of $G_1$ with $e,e_1\in S_1$.

If $e\notin M$ but $f_1\in M$, then Claim 1 implies that
$f_2\in M$.
Thus $M_1'=\{e_1\}\cup (M\cap E_1)\in \M_{1,1}$.
Clearly, $e\notin M'_1$ but $e_1\in M'_1$,
implying that
$\{e,e_1\}$ is not an equivalent set of $G_1$,
contradicting the assumption that
$S_1$ is an equivalent set of $G_1$ with $e,e_1\in S_1$.

Hence Claim 2 holds and (ii) follows.

(iii). It is trivial to verify.

(iv). Clearly, when both $G_1$ and $G_2$ are $r$-regular,
$G$ is also $r$-regular.
Assume that both $G_1$ and $G_2$ are $r$-regular graphs of class 1.
Then the edge set of each $G_i$ can be partitioned into
$r$ independent sets $E_{i,1},\cdots,E_{i,r}$.
Assume that $e_i\in E_{i,1}$ for $i=1,2$.
Then  $E(G)$ has a partition
$E_1,E_2,\cdots,E_r$
in which each subset
is an independents set of $G$, where
$$
E_1=(E_{1,1}-\{e_1\})\cup (E_{2,1}-\{e_2\})
\cup\{f_1,f_2\},
\ E_j=E_{1,j}\cup E_{2,j}, \quad \forall j=2,3,\cdots,r,
$$
implying that $G$ is of class 1.
Thus the result holds.

(v).  Suppose that $S'\ssim_G E(G)$.
Corollary~\ref{cor2-1} (ii)
implies that
$G_i-e_i-(S'\cap E_i)$ is bipartite for $i=1,2$,
a contradiction.
Thus the result holds.

(vi). Suppose that
$S'\ssim_G \emptyset$.
Corollary~\ref{cor2-1} (i)
implies that
$S'\cap E(G_i-e_i)\ssim_{G_i-e_i} \emptyset$
for $i=1,2$,
a contradiction.
Thus the result holds.
\endproof

By applying Lemma~\ref{le5-1-1}, the following conclusion
follows.

\begin{figure}[ht]
\centering
\input{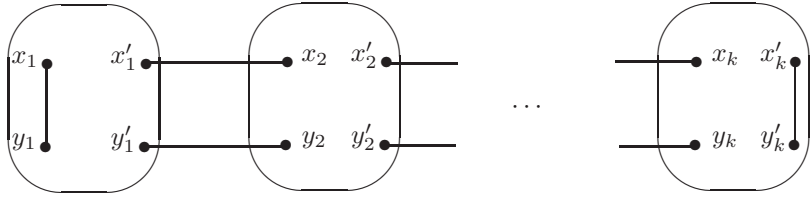}

\caption{$H_1=G_1$ and
$H_{j+1}$ is the graph
$H_{j}\#_{e'_{j},e_{j+1}} G_{j+1}$
for $j=1,2,\cdots,k-1$}
\relabel{f5}
\end{figure}

\begin{thm}\relabel{th5-1-1}
Let $G_1, G_2, \cdots, G_k$ be vertex-disjoint
$2$-connected and $r$-regular graphs of class 1
and let $S_i$ be an equivalent set of $G_i$ with $\{e_i,e'_i\}\subseteq S_i$,
where $e_i=x_iy_i$ and $e'_i=x'_iy'_i$,
for all $i=1,2,\cdots,k$.
Let $H_1=G_1$ and
let $H_{j+1}$ be the graph
$H_{j}\#_{e'_{j},e_{j+1}} G_{j+1}$
for $j=1,2,\cdots,k-1$, as shown in Figure~\ref{f5}.
Then
\begin{enumerate}
\item $H_k$ is a $2$-connected and $r$-regular graph of class 1;
\item for any subset $S$ of
$\{S_1-\{e'_1\}\}\cup \{S_k-\{e_k\}\} \cup\bigcup\limits^{k-1}_{i=2}\{S_i-\{e_i,e'_i\}\}$
with $|S|\equiv 0\modtwo$,
when $G'_i-(S\cap E(G'_i))$ is not bipartite
for some $i$ with $1\le i\le k$
and $(S\cap E(G'_j))\nssim_{G_j-e_j} \emptyset$ holds
for some $j$ with $1\le j\le k$,
$S$ is an equivalent set of $H_k$
which belongs to $\nF^*(H_k)$,
where $G'_1=G_1-\{e'_1\}$, $G'_k=G_k-\{e_k\}$
and $G'_s=G_s-\{e_s,e'_s\}$ for $2\le s\le k-1$.
\end{enumerate}
\end{thm}

\proof (i). It follows directly from Lemma \ref{le5-1-1} (iii) and (iv).

(ii).
Let $Q=\{S_1-\{e'_1\}\}\cup \{S_k-\{e_k\}\} \cup\bigcup\limits^{k-1}_{i=2}\{S_i-\{e_i,e'_i\}\}$.
Applying Lemma \ref{le5-1-1} (ii) repeatedly
shows that $Q$ is an equivalents set of $H_k$.
As $S\subseteq Q$ and $|S|\equiv 0\modtwo$,
$S\in \nF(H_k)$ holds.
As $G'_i-(S\cap E(G'_i))$ is not bipartite for some
$i$ with $1\le i\le k$,
$S\cap E(G'_i)\nssim E(G'_i)$ holds,
implying that $S\nssim_{H_k} E(H_k)$
by Corollary~\ref{cor2-1} (ii).
As $S\cap E(G'_j) \nssim_{G'_j} \emptyset$
for some $j$ with $1\le j\le k$,
Corollary~\ref{cor2-1} (ii)
implies that $S\nssim_{H_k} \emptyset$.
Hence $S\in \nF^*(H_k)$.
\endproof

By Theorem~\ref{th5-1-1}, it can be verified easily that
the graphs constructed in \cite{wei}
give a negative answer to Problem~\ref{prob1}.

\subsection{$4$-connected and $r$-regular graphs
$G$ of class 1 with $\nF^*(G)\ne \emptyset$
\relabel{sec5-2}}

In this subsection,
we construct infinitely many
$4$-connected  $r$-regular graphs $G$ of class 1
with $\nF^*(G)\ne \emptyset$,
where $r$ is an integer with $r\ge 4$.

Let $\Psi_{r}$ be the set of $4$-connected
and $r$-regular graphs  of class 1,
each of which contains an equivalent set of size $2$.
Let $Q_r$ \relabel{Qr} denote the graph
obtained from the
complete bipartite graph $K_{r,r}$
by removing two independent edges
$a_1b_1$ and $a_2b_2$
and adding two new edges $a_1a_2$ and $b_1b_2$,
where $a_1$ and $a_2$ are vertices in one partite set of $K_{r,r}$.
Observe that  $Q_r$ is a
$r$-connected and $r$-regular graph of class 1
with  an equivalent set $\{a_1a_2,b_1b_2\}$.
Thus $Q_r\in \Psi_r$.

Let $\Psi_r^*$ be the set of graphs $H\in \Psi_r$
containing an equivalent set $\{e,e'\}$
such that $H-\{e,e'\}$ is not bipartite.
From the remark in Page~\pageref{rem5-1},
it is known that $\Psi_r^*\ne \emptyset$.

For a list $L=(G_1,G_2,\cdots,G_k)$ of vertex-disjoint graphs
in $\Psi_r$, where $k\ge 3$
and $\{e_i,e'_i\}$ is an equivalent set of $G_i$
with $e_i=x_iy_i$ and $e'_i=x'_iy'_i$
for $i=1,2,\cdots,k$,
let $\C_L$ denote  the graph obtained from
$G_1,G_2,\cdots,G_k$ by
deleting edges $e_i$ and $e'_i$ and
adding new edges $f_i$ and $f'_i$ for all $i=1,2,\cdots,k$,
where $f_i=x_iy_{i+1}$, $f_i'=x'_iy'_{i+1}$,
$y_{k+1}=y_1$ and $y'_{k+1}=y'_1$.
For any $i$ with $1\le i\le k$,
assume that
$G_i-\{e_i,e'_i\}$ is not bipartite
whenever $G_i\in S_r^*$.
An example of
$\C_L$ for $k=3$ is shown in Figure~\ref{f3}.

\begin{figure}[ht]
\centering
\input{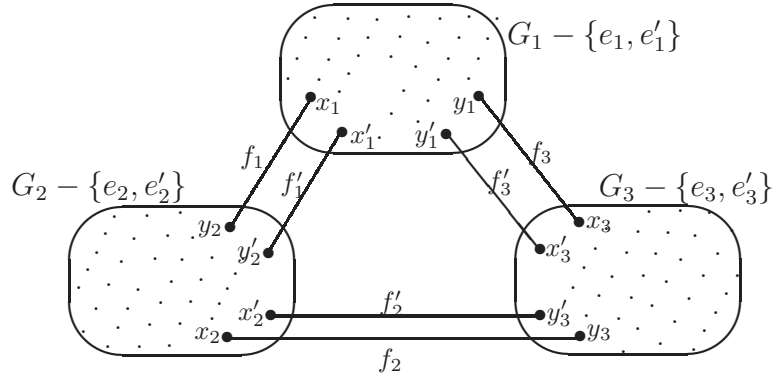}

\caption{Graph $\C_L$, where $L=(G_1, G_2,G_3)$}
\relabel{f3}
\end{figure}

\begin{lem}\relabel{le5-2-0}
Let $L=(G_1, G_2,\cdots,G_k)$ be any list of graphs
in $\Psi_r$, where $k$ is an odd number with $k\ge 3$.
The graph $\C_L$ defined above has the following properties:
\begin{enumerate}
\item $\C_L\in \Psi_r$ with equivalent sets $\{f_i,f'_i\}$
for all $i=1,2,\cdots,k$;
\item
if $G_j-\{e_j,e'_j\}$ is not bipartte
for some $j$ with $1\le j\le k$, then
$\{f_i,f'_i: i=1,2,\cdots,k\}\in \nF^*(\C_L)$ holds.
\end{enumerate}
\end{lem}

\proof
(i). As $G_i$ is $4$-connected for all
$i=1,2,\cdots,k$, it is not difficult to show that
any two non-adjacent vertices in $\C_L$
are joined by $4$ internally vertex-disjoint paths,
implying that $\C_L$ is $4$-connected.

Clearly $\C_L$ is $r$-regular.
As $G_i$ is a $r$-regular graph of class 1
and with an equivalent set $\{e_i, e'_i\}$,
$E(G_i)$ can be partitioned into
perfect matchings $E_{i,1},E_{i,2},\cdots,E_{i,r}$
with $\{e_i, e'_i\}\subseteq E_{i,1}$.
Thus, $\C_r$ is of class 1, as its edge set
can be partitioned into $r$ perfect matchings
$\E_1, \E_2,\cdots,\E_r$, where
$$
\E_1= \bigcup_{i=1}^k
\left (\{f_i, f'_i\}\cup (E_{i,1}-\{e_i, e'_i\})\right ),
\quad \E_j=\bigcup_{i=1}^k E_{i,j},\quad \forall j=2,3,\cdots,r.
$$

To show that $\{f_i,f'_i\}$ is an equivalent set of $\C_L$,
we need to apply  the following claim.

\noindent {\bf Claim 1}: For any perfect matching $M$ of $\C_L$
and any $i$ with $1\le i\le k$,
$M\cap \{f_i,f'_i\}=\{f_i\}$
implies that $M\cap \{f_{i+1},f'_{i+1}\}=\{f'_{i+1}\}$,
and $M\cap \{f_i,f'_i\}=\{f'_i\}$
implies that $M\cap \{f_{i+1},f'_{i+1}\}=\{f_{i+1}\}$.

Without loss of generality, it suffices to prove that
$M\cap \{f_1,f'_1\}=\{f_1\}$
implies $M\cap \{f_2,f'_2\}=\{f'_2\}$.
As $G_2$ is matching-covered,
$|V_2|\equiv 0\modtwo$.
Thus $M\cap \{f_1,f'_1\}=\{f_1\}$ implies that
$|M\cap \{f_2,f'_2\}|=1$.
Suppose that $M\cap \{f_2,f'_2\}=\{f_2\}$.
Then, $M_2=\{e_2\}\cup (M\cap E(G_2))$ is a perfect matching
of $G_2$.
But $e'_2\notin M_2$
contradicting the assumption
that $\{e_2,e'_2\}$ is an equivalent set of $G_2$.
Thus the claim holds.

Suppose that $\{f_i,f'_i\}$ is not an equivalent set of $\C_L$,
say $i=1$.
Then $|M\cap \{f_1,f'_1\}|=1$ holds
for some perfect matching $M$ of $\C_L$,
say $f_1\in M$ but $f'_1\notin M$.
Claim 1 implies
$M\cap \{f_2,f'_2\}=\{f'_2\}$,
$M\cap \{f_3,f'_3\}=\{f_3\}$ and so on.
As $k$ is odd, we have $M\cap \{f_k,f'_k\}=\{f_k\}$.
However, by Claim 1,
$M\cap \{f_k,f'_k\}=\{f_k\}$
implies that $M\cap \{f_1,f'_1\}=\{f'_1\}$,
a contradiction.
Hence (i) holds.

(ii).
Suppose that $G_j-\{e_j,e'_j\}$ is not bipartite for some
$j$ with $1\le j\le k$.

Let $S=\{f_i,f'_i: i=1,2,\cdots, k\}$.
As $\{f_i,f'_i\}$ is an equivalent set of $\C_L$
for all $i=1,2,\cdots,k$,
$|S\cap M|$ is even for all perfect matchings $M$ of $\C_L$,
implying that
$S\in \nF(\C_L)$ holds.

As $G_j-\{e_j,e'_j\}$ is not bipartite for some $j$ with
$1\le j\le k$,
Corollary~\ref{cor2-1} (ii) implies that
$S\nssim_{C_L} E(\C_L)$.
Suppose that $S\ssim_{\C_L} \emptyset$.
Then Proposition~\ref{pro2-1} (i) implies that
$S=\nabla_{\C_L}(U)$ for some $U\subset V(\C_L)$.
As $G_i-\{e_i,e'_i\}$ is connected,
we have $V(G_i)\subseteq U$ or
$V(G_i)\subseteq V(\C_L)-U$
for all $i=1,2,\cdots,k$.
Assume that $V(G_1)\subseteq U$.
Then  $S=\nabla_{\C_L}(U)$ implies
$V(G_2)\subseteq V(\C_L)-U$,
$V(G_3)\subseteq U$
and so on.
Since $k$ is odd,
$V(G_k)\subseteq U$,
contradicting the assumption that
$f_k,f'_k\in S=\nabla_{\C_L}(U)$.

Hence $S\in \nF^*(\C_L)$ and (ii) holds.
\endproof

By Lemma~\ref{le5-2-0}, we can prove
the following result.

\begin{cor}\relabel{cor5-2-1}
$\Psi_r^*$ is an infinite set.
\end{cor}

\proof
Let $\L$ be the family of
lists $L=(G_1,G_2,\cdots,G_k)$,
where $k\ge 3$ is odd,
$G_i\in \Psi_r$ for $i=1,2,\cdots,k$
and $G_j\in \Psi_r^*$ for at least one $j$
with $1\le j\le k$.
By the remark in Page~\pageref{rem5-1},
$\Psi_r^*\ne \emptyset$.
Thus $\L\ne \emptyset$.
By Lemma~\ref{le5-2-0},
$\C_L\in \Psi_r$ holds
for any list $L\in \L$.
Furthermore, as $G_j\in \Psi_r^*$ holds for at least one $j$,
$G_j-\{e_j,e'_j\}$ is not bipartite
for an equivalent set $\{e_j,e'_j\}$,
implying that $\C_L-\{f_i,f'_i:1\le i\le k\}$ is not
bipartite.
By Lemma~\ref{le5-2-0} (i),
$\{f_i,f'_i\}$ is an equivalent set of $\C_L$
for any $i$ with $1\le i\le k$,
implying that
$\C_L\in \Psi_r^*$ holds.
Clearly, $\C_L$ is different from anyone in the list
of $L$.
Applying Lemma~\ref{le5-2-0} repeatedly implies that
the result holds.
\endproof

By Lemma~\ref{le5-2-0} and Corollary~\ref{cor5-2-1},
we get the following result.

\begin{thm}\relabel{th5-2-1}
For any $r\ge 4$,
there are infinitely many $4$-connected and
$r$-regular graphs $H$ of class 1
with $\nF^*(H)\ne \emptyset$.
\end{thm}

\subsection{$r$-connected and $r$-regular graphs
$G$ of class 1 with $\nF^*(G)\ne \emptyset$
\relabel{sec5-3}}

For any integer $r$ with $r\ge 3$,
let $\Phi_r$ be the set of $r$-connected
and $r$-regular graphs of class 1.
Clearly, $\Phi_r$ includes
the complete bipartite graph $K_{r,r}$,
the graph $Q_r$ defined in Page~\pageref{Qr}
and the complete graph $k_{r+1}$ when $r$ is odd.

For any set $S=\{G_1, G_2,\cdots,G_r\}$ of
$r$ vertex-disjoint graphs in $\Phi_r$
with $w_i\in V(G_i)$
and $N_{G_i}(w_i)=\{v_{i,j}:  j=1,2,\cdots,r\}$
for $i=1,2,\cdots,r$,
let $\X_S$ denote the graph obtained
from $G_1-w_1,G_2-w_2,\cdots,G_r-w_r$
by adding vertices $u_1,u_2,\cdots,u_r$
and adding edges joining
$u_j$ to vertex $v_{i,j}$
for all $i=1,2,\cdots,r$ and $j=1,2,\cdots,r$,
without referring to vertices $w_i$ in $G_i$
for $i=1,2,\cdots,r$.
An example of $\X_S$ when $r=3$ is given in Figure~\ref{f4},
where $S=\{G_1,G_2, G_3\}$ and
$G_i\cong K_4$ for all $i=1,2,3$.

\begin{figure}[ht]
\centering
\input{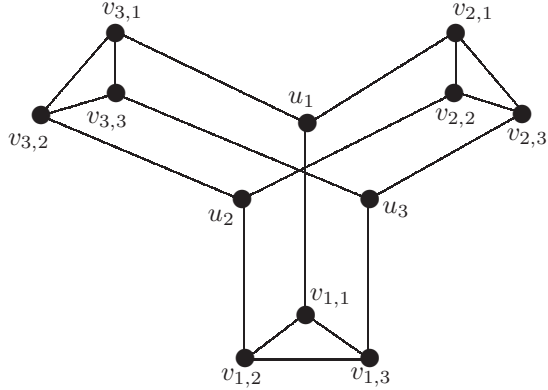}

\caption{Graph $\X_S$ with $S=\{G_1,G_2,G_3\}$
and $G_i\cong K_4$ for $i=1,2,3$}
\relabel{f4}
\end{figure}

\begin{lem}\relabel{le5-3-1}
For any set $S=\{G_1,G_2,\cdots, G_r\}$
of graphs in $\Phi_r$,
the graph $\X_S$ constructed above has the following
properties:
\begin{enumerate}
\item $\X_S\in \Phi_r$;
\item for any $i=1,2,\cdots,r$,
if both $G_i-w_i$ and $G_{j}-w_j$ are not bipartite
for some $j\in \{1,2,\cdots,r\}-\{i\}$,
then $E(G_i-w_i)\in \nF^*(\X_S)$ holds.
\end{enumerate}
\end{lem}

\proof
(i).  Observe that $\X_S$ is $r$-connected
by the two facts below:
\begin{enumerate}
\item[(a)] If both graphs $H_1$ and $H_2$ are $r$-connected
and vertex-disjoint  with
$x_i\in V(H_i)$ and  $N_{H_i}(x_i)
=\{z_{i,j}:j=1,2,\cdots,r\}$ for $i=1,2$,
then the graph obtained from $H_1-x_1$ and $H_2-x_2$
by adding edges joining $x_{1,j}$ and $x_{2,j}$ for all
$j=1,2,\cdots,r$
is also $r$-connected;

\item[(b)] for any $r$-connected graph $H$
and any $r$ independent edges $e_1,e_2,\cdots,e_r$,
the graph obtained from $H$ by
subdividing each $e_i=y_{i,1}y_{i,2}$
with a vertex, denoted by $q_i$,
adding $r-2$ new vertices
$z_1,z_2,\cdots,z_{r-2}$
and adding new edges joining
$z_j$ to $q_i$ for all $j=1,2,\cdots,r-2$
and all $i=1,2,\cdots,r$
is also $r$-connected.
\end{enumerate}
The two facts above can be verified by
proving that each pair of non-adjacent vertices are joined
by $r$ internally vertex-disjoint paths.

As $G_i$ is $r$-regular,
by the definition,
$\X_S$ is also $r$-regular.

For $i=1,2,\cdots,r$,
as $G_i$ is a $r$-regular graph of class 1,
$G_i$ has a $r$-edge-coloring which partitions
$E(G_i)$ into $r$ perfect matchings
$E_{i,1},E_{i,2},\cdots,E_{i,r}$ of $G_i$.
Assume that $w_iv_{i,j}\in E_{i,j}$ for all
$i=1,2,\cdots,r$ and $j=1,2,\cdots,r$.
Let $\pi_1,\pi_2,\cdots,\pi_r$ be permutations
of $1,2,\cdots,r$ such that
$\{\pi_s(i):s=1,2,\cdots,r\}=\{1,2,\cdots,r\}$ holds
for all $i=1,2,\cdots,r$.
Certainly such permutations exist.
Then $\E_1,\E_2,\cdots,\E_r$ defined below
form a partition of $E(\X_S)$ each of which is a matching of $\X_S$:
$$
\E_s=\bigcup_{i=1}^r \left ((E_{i,\pi_s(i)}-\{w_iv_{i,\pi_s(i)}\})
\cup \{u_iv_{i,\pi_s(i)}\}\right ),
\qquad \forall s=1,2,\cdots,r.
$$
Hence $\X_S$ is of class 1 and $\X_S\in \Phi_r$.

(ii).  For $i=1,2,\cdots,r$,
as $G_i$ is a $r$-regular graph of class 1,
$G_i$ is matching-covered,
implying that
$|V(G_i)|\equiv 0\modtwo$.
Thus $|V(G_i-w_i)|\equiv 1\modtwo$
for all $i=1,2,\cdots,r$.

For $i=1,2,\cdots,r$,
let $W_i=E(G_i-w_i)$
and $N_i=\{u_jv_{i,j}: j=1,2,\cdots,r\}$.
As $|V(G_i-w_i)|\equiv 1\modtwo$,
$|M\cap N_i|\ge 1$ holds for each perfect matching
$M$ of $\X_S$
and all $i=1,2,\cdots,r$.
But $|M\cap (N_1\cup N_2\cup\cdots \cup N_r)|=r$,
implying that
$|M\cap N_i|= 1$ holds for each perfect matching $M$ of $\X_S$
and all $i=1,2,\cdots,r$.
Thus
$|M\cap W_i|= |V(G_i)|/2-1$
holds for each perfect matching $M$ of $\X_S$,
implying that
$W_i\in \nF(\X_S)$ for all  $i=1,2,\cdots,r$.

If both $G_i-w_i$ and $G_{j}-w_j$ are not bipartite,
where $j\ne i$,
Corollary~\ref{cor2-1}
implies that $W_i\nssim_{\X_S} \emptyset$
and $W_i\nssim_{\X_S} E(\X_S)$.
Thus $W_i\in \nF^*(\X_S)$.
\endproof

For any $r\ge 3$,
let $\Phi_r^*$ be the set of graphs  $G\in \Phi_r$
such that $G-w$ is not bipartite
for every vertex $w$ in $G$.
Clearly, $Q_r\in \Phi_r^*$ and
when $r$ is odd, $K_{r+1}\in \Phi_r^*$.

\begin{lem}\relabel{le5-3-2}
For any integer $r$ with $r\ge 3$,
$\Phi_r^*$ is an infinite set.
\end{lem}

\proof Note that
$\Phi_r^*\ne \emptyset$ for any $r\ge 3$.

If $S=\{G_1, G_2,\cdots, G_r\}$ is
a set of vertex-disjoint graphs in $\Phi_r$ and $G_i\in \Phi_r^*$
holds for some pair $i,j$
with $1\le i<j\le r$,
then $\X_S-w$ is not bipartite for each vertex $w$ in $\X_S$.
By Lemma~\ref{le5-3-1},
$\X_S\in \Phi_r^*$ holds.
Note that $\X_S$ is different from any one in $S$.
Thus the result holds by applying Lemma~\ref{le5-3-1} repeatedly.
\endproof

\vspace{0.2 cm}

\noindent {\bf Remark}: \relabel{rem5-1}
By the definition in Page~\pageref{Qr},
$Q_r$ is a graph in $\Phi_r$
with an equivalent set $\{a_1a_2,b_1b_2\}$.
For any $S=\{G_1,G_2,\cdots,G_r\}$, where $G_i\in \Phi_r$
for all $i=1,2,\cdots,r$,
if $G_1$ is the graph $Q_r$ and $w_1\notin \{a_1,a_2,b_1,b_2\}$,
then it is not difficult to verify that
$\{a_1a_2,b_1b_2\}$ is an equivalent set of $\X_S$.
Furthermore, if $G_j-w_j$ is not bipartite graph
for some $j$ with $2\le j\le r$,
then $\X_S-\{a_1a_2,b_1b_2\}$
is not bipartite,
implying that $\X_S\in \Psi_r^*$ when $r\ge 4$.

\vspace{0.2 cm}

Theorem~\ref{main3}
follows directly from Lemmas~\ref{le5-3-1}
and~\ref{le5-3-2}.

\section*{Acknowledgements}
The work is partially supported by the China Scholarship Council for financial support
and  NTU AcRF project (RP 3/16 DFM) of Singapore.

\end{document}